\documentclass[12pt]{article}

\usepackage{amsmath}
\usepackage{amsfonts}
\usepackage{bigints}
\usepackage{amssymb}
\usepackage{mathabx}

\begin{document}

\renewcommand{\P}{\par \rm Proof:\ }
\newcommand{\Pe}{$\hfill{\Box}$\bigskip}
\newcommand{\Pes}{$\hfill{\Box}$}
\newcommand{\A}{\mbox{${{{\cal A}}}$}}


\author{Attila Losonczi}
\title{Dimension structures}

\date{14 December 2017}

\newtheorem{thm}{\qquad Theorem}[section]
\newtheorem{prp}[thm]{\qquad Proposition}
\newtheorem{lem}[thm]{\qquad Lemma}
\newtheorem{cor}[thm]{\qquad Corollary}
\newtheorem{rem}[thm]{\qquad Remark}
\newtheorem{ex}[thm]{\qquad Example}
\newtheorem{df}[thm]{\qquad Definition}
\newtheorem{prb}{\qquad Problem}

\maketitle

\begin{abstract}

\noindent

We are going to introduce a new algebraic, analytic structure that is a kind of generalization of the Hausdorff dimension and measure. We give many examples and study the basic properties and relations of such systems.

\noindent
\footnotetext{\noindent
AMS (2010) Subject Classifications:  28A78\\

Key Words and Phrases: Hausdorff dimension and measure}

\end{abstract}

\section{Introduction}
In this paper we introduce a new algebraic, analytic structure. In some sense it is a kind of generalization of the Hausdorff dimension and measure. We are going to create the basic building blocks of this new theory while presenting many examples which had inspired the development of such area.

Roughly speaking a dimension structure is a system that consists of a set and many measurements which measure the points of the set. However most of those measurements measure a point to either $0$ or $+\infty$. Moreover the measurements are (partially) ordered and there is a rule which says roughly that the greater the measurement position the smaller its value. More precisely if the value is not infinity then all greater measurements measure the point to 0. Which gives the opportunity to define dimension of each point. Hence the points can be compared by two values, the dimension and the measurement at the dimension.

We present many examples for various dimension structures. Actually there are example groups like geometric measure theory types, accumulation point types, astronomical distance types, convergence speed types and restricted Lebesgue measure types.

In the last section we deal with the task how one can build new dimension structures from already existing ones. E.g. we define substructure, quotient structure of dimension structures and also sum and direct product of dimension structures and the more exotic measure sum. Then we devote a subsection to mappings between dimension structures. Finally we investigate extensions i.e. how a pre-dimension structure can be extended to a dimension structure.

\subsection{Basic notions and notations}\label{ssbn}

We enumerate many notions and notations that we will frequently use later.

\smallskip

For a set $X$ let $P(X)$ denote the set of all subsets of $X$.

\smallskip

Let $\overline{\mathbb{R}^+}=\mathbb{R}^+\cup\{0,+\infty\}=[0,+\infty]$.

\smallskip

Throughout this subsection (\ref{ssbn}) let $S$ denote a partially ordered set (i.e. poset).

If a partially ordered set $S$ has a minimum (maximum) element we also denote it by $-\infty$ ($+\infty$). If the set does not have a minimum (maximum) element then the symbol $-\infty$ ($+\infty$) denotes a new element that is less than (greater than) all elements in $S$.
Let $\bar S=S\cup\{-\infty,+\infty\}$. E.g. if $S$ has a minimum and maximum then $S=\bar S$.

\smallskip

Within a partially ordered set $S$ let us apply the usual convension that $\inf\emptyset=+\infty,\ \sup\emptyset=-\infty$.

\medskip
Let us recall some usual definitions regarding posets.

$S$ is called \textbf{ordered} if $\forall x,y\in S$ are comparable that is $x\leq y$ or $y\leq x$ holds.

\smallskip

$S$ is called \textbf{dense} if $\forall x<y\in S\ \exists z\in S$ such that $x<z<y$.

\smallskip

$S$ is called \textbf{discrete} if the topology induced by the partial order is discrete. Equivalently if there is no infinite chain between two comparable points.

\smallskip

$S$ is called \textbf{complete} if $\forall P\subset S\ \sup P$ and $\inf P$ exist.

\smallskip

$S$ is a \textbf{lattice} if $\forall s,p\in S\ \inf\{s,p\}$ and $\sup\{s,p\}$ exist.

\smallskip

$S'\subset S$ is called a principal filter if $S'=\{s\in S:\inf S'\leq s\}$. It implies that $\inf S'=\min S'\in S$.

\smallskip

If $s\in S$ has a successor then it is denoted by $s^+=\min\{p\in S:p>s\}$ while the predecessor (if exists) by $s^-=\max\{p\in S:p<s\}$.

\smallskip

$P\subset S$ is \textbf{convex} if $a,c\in P,b\in S,a<b<c$ implies that $b\in P$. 

$P\subset S$ is called \textbf{up-convex} if $a\in P,b\in S,a<b$ implies that $b\in P$. An up-convex set is convex. If $+\infty\in P$ and $P$ is convex then it is up-convex.

$P\subset S$ is called \textbf{down-convex} if $a\in P,b\in S,b<a$ implies that $b\in P$. A down-convex set is convex. If $-\infty\in P$ and $P$ is convex then it is down-convex.

Obviously if $\pm\infty\in S$ then the only up- and down-convex set is $S$.

\medskip

Let us define some usual operations with $+\infty$:
$$0\cdot(+\infty)=0,\ (+\infty)+(+\infty)=+\infty,\ (+\infty)\cdot(+\infty)=+\infty.$$
If $r\in\mathbb{R}^+$ then $$r\cdot (+\infty)=+\infty,\ \frac{+\infty}{r}=+\infty,\ r+(+\infty)=+\infty.$$

\smallskip

Let us extend the domain of the $sign$ function with $+\infty$:
$$sign(x)=
\begin{cases}
-1&\text{if }x<0\\
0&\text{if }x=0\\
1&\text{if }x>0\\
+\infty&\text{if }x=+\infty.
\end{cases}$$

\smallskip

In this paper $\lambda$  and ${\cal{L}}$ denote the Lebesgue measure on $\mathbb{R}$ and the set of Lebesgue measurable sets respectively.
\section{The concept of dimension structure}
\begin{df}${\cal{D}}=\langle X,\mu_s:s\in S\rangle$ is a dimension structure if $X$ is a set, $\langle S,<\rangle$ is a partially ordered set, $\mu_s\ (s\in S)$ is a function $\mu_s:X\to\mathbb{R}^+\cup\{0,+\infty\}$ and the following axioms hold:
\begin{itemize}
\item[(ax1)] If $x\in X,s,p\in S,s<p$ and $\mu_s(x)<+\infty$ then $\mu_p(x)=0$.
\item[(ax2)] If $x\in X,s,p\in S$,  $0<\mu_s(x)<+\infty$ and $\mu_p(x)<+\infty$ then $s,p$ are comparable elements.
\item[(ax3)] $\forall x\in X\ \exists\inf\{s\in S:\mu_s(x)=0\}$.
\end{itemize}

If $-\infty\not\in S$ then set $\mu_{-\infty}(x)=+\infty\ (x\in X)$. When $+\infty\not\in S$ then set $\mu_{+\infty}(x)=0\ (x\in X)$.
\Pes
\end{df}

\begin{rem}With the extension mentioned in the definition $\mu$ can be considered as a function $\mu:\bar{S}\times X\to[0,+\infty]$.
\end{rem}

\begin{rem}Clearly (ax2) is statisfied if $\langle S,<\rangle$ is ordered.

Obviously (ax3) holds if $\langle S,<\rangle$ is a complete partially ordered set.

When $\langle S,<\rangle$ is a complete ordered set then only (ax1) has to be checked.
\end{rem}

\begin{rem}Examining (ax3) it can happen that $\inf\{s\in S:\mu_s(x)=0\}\not\in S$ because it equals to $-\infty\in \bar{S}-S$. Or similarly when $\{s\in S:\mu_s(x)=0\}$ is empty then the infimum is $+\infty$ that may not be in $S$.
\end{rem}

\begin{df}Let $\langle X,\mu_s:s\in S\rangle$ be a dimension structure and $x\in X$. We will use the following notations: $$S_x=\{s\in S:\mu_s(x)<+\infty\},\ S^0_x=\{s\in S:\mu_s(x)=0\},\ S^{+\infty}_x=S-S_x.$$
\end{df}

With that (ax3) can be formulated in a way that infimum of $S^0_x$ always exists.

\begin{df}Let $\langle X,\mu_s:s\in S\rangle$ be a dimension structure and $x\in X.$ If $0<\mu_s(x)<+\infty$ for some $s\in S$ then $x$ is called an $s$-point. 

Let us call $x$ a dimension-point or dim-point if it is an $s$-point for some $s\in S$.
\end{df}

In the sequel we are going to enumerate facts that are the analog to ones in the theory of Hasudorff dimension and measure. In the following statements we assume that a dimension structure is given with the above attributes.

\begin{prp}\label{p0i}If $x\in X,\ s<p\ (s,p\in S),\ 0<\mu_p(x)$ then $\mu_s(x)=+\infty$.
\end{prp}
\P Indirect, using (ax1).
\Pes

\begin{cor}\label{cii}If $x\in X,\ s<p\ (s,p\in S),\ \mu_p(x)=+\infty$ then $\mu_s(x)=+\infty$. Equivalently $S_x^{+\infty}$ is down-convex $(x\in X)$. \Pes
\end{cor}

\begin{prp}\label{pzi}If $x\in X,\ s<p\ (s,p\in S),\ \mu_s(x)=0$ then $\mu_p(x)=0$. Equivalently $S_x^0$ is up-convex $(x\in X)$.
\end{prp}
\P Apply (ax1).
\Pe


The following statements prepare the concept of dimension.

\begin{cor}\label{pddsd}For a given $x\in X$ there can be at most one $s\in S$ such that $0<\mu_s(x)<+\infty$. 
\end{cor}
\P By (ax2) if there were two then they have to be comparable and then (ax1) yields a contradiction.
\Pes

\begin{cor}(ax3) is equivalent to $\forall x\in X\ \exists\inf S_x$.
\end{cor}

\begin{cor}\label{csic}For a given $x\in X$ 
$$\sup\{s\in S:\mu_s(x)=+\infty\}\leq\inf\{s\in S:\mu_s(x)=0\}$$
if the $\sup$ and $\inf$ exist and they are comparable. 
\end{cor}
\P Let $s=\sup\{s\in S:\mu_s(x)=+\infty\},\ i=\inf\{s\in S:\mu_s(x)=0\}$. Let suppose that $i<s$. Then there is $a\in X,\ i<a\leq s$ such that $\mu_a(x)=+\infty$. And similarly there is a $b\in X,\ i\leq b< a$ such that $\mu_b(x)=0$. But then by \ref{pzi} $\mu_a(x)=0$ which is a contradiction.
\Pes

\begin{cor}Let ${\cal{D}}=\langle X,\mu_s:s\in S\rangle$ be a dimension structure with $S$ being ordered. Then for a given $x\in X$ if $\mu_s(x)=+\infty$ and $\mu_p(x)=0$ then $s<p$. Moreover if $0<\mu_q(x)<+\infty$ then $s<q<p$. \Pes
\end{cor}

\begin{prp}Let ${\cal{D}}=\langle X,\mu_s:s\in S\rangle$ be a dimension structure with $S$ being dense, complete and ordered. Then 
$$\sup\{s\in S:\mu_s(x)=+\infty\}=\inf\{s\in S:\mu_s(x)=0\}.$$
\end{prp}

\P If $S_x^{+\infty}=\emptyset$ then by \ref{p0i} $\inf S_x^0=-\infty$ hence the assertion holds. The case $S_x^0=\emptyset$ is similar.

Suppose that $S_x^{+\infty}\ne\emptyset$ and $\ne S$. Then by \ref{csic} $\sup S_x^{+\infty}\leq\inf S_x^0$ holds. If they would not be equal then there were infinitely many points between them which is a contradiction since all but one point satisfy $\mu_s(x)=0$ or $+\infty$.  
\Pes

\begin{df}Let ${\cal{D}}=\langle X,\mu_s:s\in S\rangle$ be a dimension structure. For $x\in X$ set  
$$\dim_{\cal{D}} x=\dim x=\inf\{s\in S:\mu_s(x)<+\infty\}\in \bar S.$$
It is called the dimension of $x$ in ${\cal{D}}$.
\end{df}

\begin{rem}\label{rdim+inf}If ${\cal{D}}=\langle X,\mu_s:s\in S\rangle$ is a dimension structure with $S$ being dense, complete and ordered then
$$\dim_{\cal{D}} x=\sup\{s\in S:\mu_s(x)=+\infty\}\ (x\in X).$$
\end{rem}

In a dimension structure we can easily define a partial order on $X$ using dimension and the given measure.

\begin{df}Let $\langle X,\mu_s:s\in S\rangle$ be a dimension structure. We define a partial order on $X$. Let $x,y\in X$.
$$x\leq_{{\cal{D}}} y\iff
\begin{cases}
\dim x<\dim y \text{ or } \\
\dim x=\dim y=d \text{ and } \mu_d(x)\leq\mu_d(y).
\end{cases}$$
Obviously if $S$ is ordered then it is an order too.
\end{df}

If we want to measure the size of an object $x\in X$,  $\mu_{\dim x}(x)$ is not enough alone, actually we need $\dim x$ as well. Hence we can express the size of $x$ by two numbers $\dim x$ and $\mu_{\dim x}(x)$.

\begin{df}\label{dmu}Let ${\cal{D}}=\langle X,\mu_s:s\in S\rangle$ be a dimensional structure. For $x\in X$ let $\mu_{\cal{D}}:X\to S\times[0,+\infty]$ $$\mu_{\cal{D}}(x)=(\dim x,\mu_{\dim x}(x)).$$ 
\end{df}

\begin{prp}Let ${\cal{D}}=\langle X,\mu_s:s\in S\rangle$ be a dimensional structure. If we equip $S\times[0,+\infty]$ with the lexicographic order then 
$$x\leq_{{\cal{D}}} y\iff \mu_{\cal{D}}(x)\leq\mu_{\cal{D}}(y)\ (x,y\in X).$$ \Pes
\end{prp}

\begin{df}For dimension structure ${\cal{D}}$ and $d\in S, m\in[0,+\infty]$ define $C^{\cal{D}}_{d,m}=C_{d,m}=\{x\in X:\dim x=d,\mu_d(x)=m\}$.
\end{df}

\begin{prp}$x\in X$ is an $s$-point iff $\dim x=s$ and $0<\mu_{\dim x}(x)<+\infty$. \Pes
\end{prp}

\subsection{Special classes of dimension structures}

Here we define some important properties of dimension structures.

\begin{df}
We call a dimension structure ${\cal{D}}=\langle X,\mu_s:s\in S\rangle$ fully-normal if $\forall s\in S\ \exists x\in X$ such that $0<\mu_s(x)<+\infty$ i.e. $x$ is an $s$-point. ${\cal{D}}$ is normal if $\forall s\in S,s\notin\{-\infty,+\infty\}$ there exists an $s$-point. ${\cal{D}}$ is called quasi-normal if $\forall s\in S,s\notin\{-\infty,+\infty\}$ there exists $x\in X$ such that $\dim x=s$.
\end{df}

\begin{df}
We call a dimension structure $\langle X,\mu_s:s\in S\rangle$

p-strong if $\forall x\in X\ \exists s\in S$ such that $\mu_s(x)<+\infty$,

m-strong if $\forall x\in X\ \exists s\in S$ such that $0<\mu_s(x)$.

If it is p-strong and m-strong at the same time then it is called strong.
\end{df}

\begin{prp}Propery p-strong is equivalent to the condition that $\dim x\ne+\infty\ (x\in X)$, while property m-strong is equivalent to the condition that $\dim x\ne-\infty\ (x\in X)$.
\end{prp}

We might have added the following property to the axioms because it is so generic however it might have restricted the concept of dimension structure too much. 

\begin{df}A dimension structure ${\cal{D}}=\langle X,\mu_s:s\in S\rangle$ is called principal if $\forall x\in X\ S_x\cup\{\inf S_x\}$ is a principal filter in $S$. 
\end{df}

\begin{rem}${\cal{D}}$ is principal iff $S_x-\{\inf S_x\}=\{s\in S:s>\inf S_x\}$.
\end{rem}

\begin{rem}If $S$ is ordered then ${\cal{D}}$ is principal.
\end{rem}

\begin{prp}\label{pprinceq}Dimension structure ${\cal{D}}=\langle X,\mu_s:s\in S\rangle$ is principal iff $x\in X,\ s>\dim x$ implies that $\mu_s(x)=0$.
\end{prp}
\P Necessity is obvious. We prove sufficiency. Let us assume that ${\cal{D}}$ is not principal. Then there is $x\in X,s\in S,s>\dim x$ such that $\mu_s(x)=+\infty$ which is a contradiction.
\Pes

\begin{ex}If ${\cal{D}}$ is not principal then it can happen that for some $x\in X$ there is $s\in S$ such that $s>\dim x$ and $\mu_s(x)=+\infty$.

To provide such example let $X=\{x_0\},\ S=\{-n:n\in\mathbb{N}\}\cup\{-\infty,a\}$ where $a$ is a new symbol differing from the other ones. Set 
$$s<p\iff
\begin{cases}
-s,-p\in\mathbb{N}\text{ and } s<p\\
s=-\infty\\
s=a,p=-1.
\end{cases}$$
Then $S$ is a lattice. Set
$$\mu_s(x_0)=
\begin{cases}
+\infty&\text{if }s\in\{-\infty,a\}\\
0&\text{otherwise}.
\end{cases}$$
Then $\langle X,\mu_s:s\in S\rangle$ is a dimensional structure and $\dim x_0=-\infty$ and obviously the assertion holds since $\mu_a(x_0)=+\infty$. 
\end{ex}

In a plain dimension structure among the points of $X$ there is no any connection. The following definition provides some.

\begin{df}Let ${\cal{D}}=\langle X,\mu_s:s\in S\rangle$ be a dimensional structure and let $\langle X,\leq\rangle$ be partially ordered too. If $\alpha$ is a cardinal number, we say that ${\cal{D}}$ is $\alpha$-synchronized if the followings hold:
\begin{itemize}
\item[1.] If $x\leq y\ (x,y\in X)$ then $\mu_{\cal{D}}(x)\leq\mu_{\cal{D}}(y)$.

\item[2.] If $Y\subset X,\ |Y|\leq\alpha$ and $\sup Y$ exists in $X$ then $\sup\{\dim y:y\in Y\}$ exists as well and $\dim(\sup Y)=\sup\{\dim y:y\in Y\}$ holds.
\end{itemize}
We say that ${\cal{D}}$ is finitely-synchronized if condition 2 holds for finite sets $Y$ (and condition 1 remains valid too).
\end{df}

\begin{rem}Alone condition 1 gives that if $Y\subset X$, $\sup Y$ and $\sup\{\dim y:y\in Y\}$ exist then $\sup\{\dim y:y\in Y\}\leq\dim(\sup Y)$.
\end{rem}
\P Condition 1 gives that $x\leq y$ implies that $\dim x\leq\dim y$. And if $y\in Y$ then $y\leq\sup Y$ hence $\dim y\leq\dim\sup Y$.
\Pes

\begin{prp}\label{psync1}Let ${\cal{D}}=\langle X,\mu_s:s\in S\rangle$ be a dimensional structure with $S$ being ordered. Then condition 1 is equivalent to 
\begin{equation}\label{eq1}\text{if }x\leq y\text{ then }\forall s\in S\ \mu_s(x)\leq\mu_s(y).\end{equation}
\end{prp}
\P Let us prove necessity. If $\dim x=\dim y$ then it is obvious. Let $\dim x<\dim y$ and assume that $\exists s\in S$ such that $\mu_s(x)>\mu_s(y)$. It gives that $\mu_s(y)<+\infty$ and then $\dim y\leq s$. But it implies that $\dim x< s$ which yields that $\mu_s(x)=0$ that is a contradiction because $\mu_s(y)$ should be negative.

Let us show sufficiency. The new condition (\ref{eq1}) gives that $$\{s\in S:\mu_s(y)<+\infty\}\subset\{s\in S:\mu_s(x)<+\infty\}$$ which implies that $\dim x\leq\dim y$. If $\dim x<\dim y$ the we are done, while if $\dim x=\dim y=d$ then $\mu_d(x)\leq\mu_d(y)$ gives the statement.
\Pes

\begin{rem}There is a weaker form of Condition1:

1'. If $x\leq y\ (x,y\in X)$ then $\dim x\leq\dim y$.

\medskip

If $\langle X,\mu_s:s\in S\rangle$ is a dimensional structure with $X$ being ordered then for finite synchronization Condition 1' is sufficient.
\end{rem}

\begin{prp}\label{psynf}Let ${\cal{D}}=\langle X,\mu_s:s\in S\rangle$ be a dimensional structure with $X$ being ordered. Then condition 1 implies finite synchronization. \Pes
\end{prp}

The question arises that if we equipped $X$ with the partial order from $\leq_{\cal{D}}$ could this make ${\cal{D}}$ (at least) finitely synchronized. The answer is negative even when $S$ and $\langle X,\leq_{\cal{D}}\rangle$ are lattices as the next example shows.

\begin{ex}Let $S=\{a,b,c,d,e\},X=\{x,y,z,w\}$. The order on $S$: $a<c<d<e,a<b<d$. Clearly $S$ is a lattice. Let us define $\mu$: $\forall s\in S\ \mu_s(x)=0,$ if $s\in\{c,d,e\}$ then $\mu_s(z)=0$, if $s\in\{b,d,e\}$ then $\mu_s(y)=0$ and $\mu_e(w)=0$. Let $\mu=+\infty$ for all not yet defined cases.

Then $\dim x=a,\dim y=b,\dim z=c,\dim w=e$. We get that $x\leq_{\cal{D}}y\leq_{\cal{D}}w,x\leq_{\cal{D}}z\leq_{\cal{D}}w$. Therefore $X$ becomes a lattice but $d=\sup\{\dim y,\dim z\}<\dim \sup\{y,z\}=e$.
\end{ex}

We define one more property.

\begin{df}
We call a dimension structure $\langle X,\mu_s:s\in S\rangle$

p-small if $\forall x\in X\ \mu_{\dim x}(x)<+\infty$,

m-small if $\forall x\in X\ 0<\mu_{\dim x}(x)$.

If it is p-small and m-small at the same time then it is called small and that is eqivalent to that all $x\in X$ is a $\dim x$-set.
\end{df}

\begin{prp}Let ${\cal{D}}$ be a dimension structure. If ${\cal{D}}$ is p-small then it is p-strong. If ${\cal{D}}$ is m-small then it is m-strong. If ${\cal{D}}$ is small then it is strong. \Pes
\end{prp}

\subsection{Discrete ordered dimension structures}

We mention an important special case here.

\begin{prp}Let ${\cal{D}}=\langle X,\mu_s:s\in S\rangle$ be given such that $X$ is a set, $\langle S,<\rangle$ is a partially ordered set, $\mu_s\ (s\in S)$ is a function $\mu_s:X\to\mathbb{R}^+\cup\{0,+\infty\}$ and let $S$ be discrete and ordered. Then ${\cal{D}}$ is a dimension structure iff the following holds:

(ax1') $\mu_s(x)<+\infty\Rightarrow\mu_{s^+}(x)=0\ (\text{whenever }s\in S,x\in X$ and $s^+$ denotes the successor of $s)$
\end{prp}
\P (ax1) can be shown by induction from (ax1'). (ax2) is satisfied because $S$ is ordered, and (ax3) holds because a discrete ordered set is complete.
\Pes

\begin{prp}If ${\cal{D}}=\langle X,\mu_s:s\in S\rangle$ is a dimension structure with $S$ being discrete and ordered then
$$(\sup S_x^{+\infty})^+\leq\dim x\leq\inf S_x^0.$$
\end{prp}
\P $\dim x=\min\{s\in S:\mu_s(x)<+\infty\}$
\Pes

\begin{cor}If ${\cal{D}}=\langle X,\mu_s:s\in S\rangle$ is a p-strong dimension structure with $S$ being discrete and ordered, $x\in X$ then
$\mu_{\dim x}(x)<+\infty$ i.e. ${\cal{D}}$ is p-small. \Pes
\end{cor}

\section{Examples}

\begin{ex}\label{eH}The classic Hausdorff dimension and measure on $\mathbb{R}$ is a dimension structure: ${\cal{D}}=\langle P(\mathbb{R}):\mu_s:s\in[0,1]\rangle$ where $\mu_s$ is the $s$-dimensional Hausdorff measure. Moreover if we equip $P(\mathbb{R})$ with the partial order from set inclusion then ${\cal{D}}$ becomes an $\aleph_0$-synchronized dimension structure.
\end{ex}

\begin{rem}We can also get dimension structures if we replace "Hausdorff dimension" to "box-counting dimension" or "packing dimension" or any other such geometric type dimension and measure.
\end{rem}

\begin{ex}Let us use the notation $H'$ for the accumulation points of a subset $H$ of $\mathbb{R}$. Let $H^{(0)}=H$ and $H^{(n+1)}=(H^{(n)})'$ for $n\in\mathbb{N}\cup\{0\}$. If $n\in\mathbb{N}\cup\{0\}$ let 
$$\mu_n(H)=
\begin{cases}
|H^{(n)}|&\text{if }|H^{(n)}| \text{ is finite}\\
+\infty&\text{otherwise}
\end{cases}$$

Then ${\cal{D}}=\langle P(\mathbb{R}),\mu_n:n\in\mathbb{N}\cup\{0\}\rangle$  is a dimension structure. If we equip $P(\mathbb{R})$ with the partial order from set inclusion then ${\cal{D}}$ becomes a finitely-synchronized normal dimension structure.
\end{ex}
\P Clearly if $|H^{(n)}|$ is finite then $H^{(n+1)}=\emptyset$ and this yields (ax1').

We show that ${\cal{D}}$ is finitely-synchronized. Let us recall the facts that $A\subset B\Rightarrow A'\subset B'$ and $(A\cup B)'=A'\cup B'$. These imply that $A\subset B\Rightarrow A^{(n)}\subset B^{(n)}$ and $(A\cup B)^{(n)}=A^{(n)}\cup B^{(n)}\ (n\in\mathbb{N})$ which simply give the statement.

Normality is obvious since if $\dim H\not\in \{0,+\infty\}$ then $H$ is a $\dim H$-point.
\Pes

\begin{rem}Obviously ${\cal{D}}$ is not $\aleph_0$-synchronized as the following sets testify: $H_i=\{\frac{1}{i}\}\ (i\in\mathbb{N})$ since $\dim H_i=0$ while $\dim\bigcup_1^{+\infty}H_i=1$.
\end{rem}

\begin{rem}If we slighly modify ${\cal{D}}$ with leaving out the empty set from the underlying set then $\langle P(\mathbb{R})-\{\emptyset\},\mu_n:n\in\mathbb{N}\cup\{0\}\rangle$ becomes m-small.
\end{rem}

\begin{ex}We can generalize the previous example further. Let $\alpha<\omega_1$ be an ordinal number. If $\alpha=\beta+1$ is a successor ordinal then let $H^{(\alpha)}=\big(H^{(\beta)}\big)'$, and when $\alpha$ is a limit ordinal then let $H^{(\alpha)}=\bigcap\limits_{\beta<\alpha}H^{(\beta)}$.
 If $\alpha<\omega_1$ let 
$$\mu_{\alpha}(H)=
\begin{cases}
|H^{(\alpha)}|&\text{if }|H^{(\alpha)}| \text{ is finite}\\
+\infty&\text{otherwise.}
\end{cases}$$
Then ${\cal{D}}=\langle P(\mathbb{R}),\mu_{\alpha}:\alpha<\omega_1\rangle$  is a dimension structure. If we equip $P(\mathbb{R})$ with the partial order from set inclusion then ${\cal{D}}$ becomes a finitely-synchronized normal dimension structure.

In this case the dimension of a set can be greater then $\omega_0$.
\end{ex}

\begin{ex}\label{eseqsi}Let $X=\{(x_n) \text{ sequence on } \mathbb{R}:x_n\to+\infty\}$. We want to measure how fast the sequence tends to infinity. For $\alpha>0$ set 
$$l_{\alpha}(x_n)=\varliminf\limits_{n\to\infty}\frac{x_n}{n^{\alpha}},\ L_{\alpha}(x_n)=\varlimsup\limits_{n\to\infty}\frac{x_n}{n^{\alpha}}.$$
Then $\langle X,\ l_{\alpha}:\alpha>0\rangle$ and $\langle X,\ L_{\alpha}:\alpha>0\rangle$ are normal dimension structures.
\end{ex}
\P We show it for $\langle X,\ l_{\alpha}:\alpha>0\rangle$, the other is similar.

Suppose $l_{\alpha}(x_n)<+\infty, \alpha<\beta\in\mathbb{R}$. Then $\varliminf\frac{x_n}{n^{\beta}}=\varliminf\frac{x_n}{n^{\alpha}}\frac{1}{n^{\beta-\alpha}}=0$ because $\varliminf\frac{1}{n^{\beta-\alpha}}=0$.

Obviously $\forall\alpha\ \exists (x_n)$ such that $0<l_{\alpha}(x_n)<+\infty$, e.g. $n^{\alpha}$.
\Pe

We can generalize the previous example further. (The proof is similar hence omitted.) 

\begin{ex}\label{eseqs}Let $X=\{(x_n) \text{ sequence on } \mathbb{R}:x_n\to+\infty\}$. Set $S=\omega_0\times\big(\mathbb{R}^+\cup\{0\}\big)$ and take the lexicographic order on $S$. For $m\in\omega_0$ we define functions on $\mathbb{R}$: $f_0(x)\equiv 1, f_1(x)=e^x, f_{m+1}(x)=e^{f_m(x)}$ 
If $(x_n)\in X,\ (m,\alpha)\in S$ let
$$l_{(m,\alpha)}(x_n)=\varliminf\limits_{n\to\infty}\frac{x_n}{f_m(n)\cdot n^{\alpha}},\ L_{(m,\alpha)}(x_n)=\varlimsup\limits_{n\to\infty}\frac{x_n}{f_m(n)\cdot n^{\alpha}}.$$
Then $\langle X,\ l_{(m,\alpha)}:(m,\alpha)\in S\rangle$ and $\langle X,\ L_{(m,\alpha)}:(m,\alpha)\in S\rangle$ are normal dimension structures. Now the dimension of a point is of the form $(m,\alpha)\in S$.

Of course further such generalizations are still possible.
\end{ex}

We now present a different type of generalization.

\begin{ex}\label{eseqso}Let $X=\{(x_n) \text{ sequence on } \mathbb{R}:x_n\to+\infty\}$. For $\alpha,\beta>0$ let 
$$l_{\alpha,\beta}(x_n)=
\begin{cases}
+\infty&\text{if }\varliminf\limits_{n\to\infty}\frac{x_{2n}}{(2n)^{\alpha}}=+\infty\text{ or }\varlimsup\limits_{n\to\infty}\frac{x_{2n+1}}{(2n+1)^{\beta}}=+\infty\\
\varliminf\limits_{n\to\infty}\frac{x_{2n}}{(2n)^{\alpha}}\cdot\varlimsup\limits_{n\to\infty}\frac{x_{2n+1}}{(2n+1)^{\beta}}&\text{otherwise}.
\end{cases}$$
Then $\langle X,\ l_{\alpha,\beta}:\alpha,\beta>0\rangle$ is a dimension structure where we take the product order on $\mathbb{R}^+\times\mathbb{R}^+$ (the underlying poset is not ordered any more).
\end{ex}

\P Let us check (ax1). Suppose $l_{\alpha,\beta}(x_n)<+\infty$. Then $$\varliminf\limits_{n\to\infty}\frac{x_{2n}}{(2n)^{\alpha}}<+\infty \text{ and } \varlimsup\limits_{n\to\infty}\frac{x_{2n+1}}{(2n+1)^{\beta}}<+\infty.$$ 

Let $(\alpha,\beta)<(\alpha',\beta')$. Then $\alpha<\alpha'$ and $\beta\leq\beta'$ or $\beta<\beta'$ and $\alpha\leq\alpha'$. Assume the first holds (the other is similar). 
Then $$\varliminf\limits_{n\to\infty}\frac{x_{2n}}{(2n)^{\alpha'}}=\varliminf\limits_{n\to\infty}\frac{x_{2n}}{(2n)^{\alpha}}\frac{1}{(2n)^{\alpha'-\alpha}}=0$$ because $\varliminf\frac{1}{(2n)^{\alpha'-\alpha}}=0$. Obviuosly $\varlimsup\limits_{n\to\infty}\frac{x_{2n+1}}{(2n+1)^{\beta'}}\leq\varlimsup\limits_{n\to\infty}\frac{x_{2n+1}}{(2n+1)^{\beta}}$ hence $\varlimsup\limits_{n\to\infty}\frac{x_{2n+1}}{(2n+1)^{\beta'}}<+\infty$. Which implies that $l_{\alpha',\beta'}(x_n)=0$.

\smallskip

We show that (ax2) holds too. Let $0<l_{\alpha,\beta}(x_n)<+\infty$ and $l_{\alpha',\beta'}(x_n)<+\infty$. Then $$0<\varliminf\limits_{n\to\infty}\frac{x_{2n}}{(2n)^{\alpha}}<+\infty \text{ and } 0<\varlimsup\limits_{n\to\infty}\frac{x_{2n+1}}{(2n+1)^{\beta}}<+\infty$$
and
$$\varliminf\limits_{n\to\infty}\frac{x_{2n}}{(2n)^{\alpha'}}<+\infty \text{ and } \varlimsup\limits_{n\to\infty}\frac{x_{2n+1}}{(2n+1)^{\beta'}}<+\infty.$$
From that we get that $\alpha\leq\alpha'$ and $\beta\leq\beta'$ which yields that $(\alpha,\beta)\leq(\alpha',\beta')$.

\smallskip

Let us validate (ax3). Evidently $$\{(\alpha,\beta):l_{\alpha,\beta}(x_n)<+\infty\}=\{\alpha:\varlimsup\frac{x_{2n}}{(2n)^{\alpha}}<+\infty\}\times\{\beta:\varlimsup\frac{x_{2n+1}}{(2n+1)^{\beta}}<+\infty\}$$ hence the infimum exists.
\Pes

\begin{rem}The dimension structure is not principal in \ref{eseqso}.

To show that take the sequence $x_n=n\log n$. Clearly $$\varliminf\limits_{n\to\infty}\frac{x_{2n}}{2n}=\varlimsup\limits_{n\to\infty}\frac{x_{2n+1}}{2n+1}=+\infty.$$ If $\alpha>1$ and $\beta>1$ then $$\varliminf\limits_{n\to\infty}\frac{x_{2n}}{(2n)^{\alpha}}=\varlimsup\limits_{n\to\infty}\frac{x_{2n+1}}{(2n+1)^{\beta}}=0$$ which gives that $\dim (x_n)=(1,1)$. If we take e.g. $(1,5)$ then $(1,1)<(1,5)$ and $l_{1,5}(x_n)=+\infty$ because the first factor is $+\infty$.
\end{rem}

\begin{rem}In \ref{eseqso} we can get a different however very similar example if in the definition of $l_{\alpha,\beta}(x_n)$ we change the product to $\min$.
\end{rem}

\begin{rem}We can slightly extend $X$ in the previous examples (\ref{eseqsi},\ref{eseqs},\ref{eseqso}) by adding generalized sequences which can also take the value $+\infty$.
\end{rem}

\begin{rem}Using the previous methods one can measure convergence speed as well. Let $a_n\to a\in\mathbb{R}$. Let 
$$x_n=
\begin{cases}
\frac{1}{|a_n-a|}&\text{if }a_n\ne a\\
+\infty&\text{otherwise}.
\end{cases}$$
Then $(x_n)$ is a generalized sequence tending to infinity and we can apply the generalized form of the previously definied dimensions and measures of $(x_n)$ to $(a_n)$.
\end{rem}

\begin{ex}We define a mathematical model for describing astronomical/atomic distances. First let us explain for what exactly we want to construct a model.

When we measure distances in the universe we find that there are distance classes in the following sense. E.g. we have one class of distances when we measure man size objects on Earth. Another class is when we measure distance between planetary objects (planets, meteorites, etc.) since those distances are many orders of magnitude greater than the previous ones. Then a higher class is when we measure distance between stars because again those distances are far bigger than the previous ones. Next class is the distance between galaxys.  

Or we can go downward as well. When we measure distances between components in a cell than those are many order of magnitude smaller than distances of man size objects. The next class down can be the distance between atoms, then distance between elementary particles.

And each time in a class we get distances that do not differ from each other much, while distances in different classes are differ by many orders of magnitude.

Let us now present a mathematical model for such system.

\medskip

Let $X$ be a set. For each $n\in\mathbb{Z}$ let $\rho_n$ be a pseudo--metric on $X$ with the additional property that $\rho_n(x,y)$ can take the value $+\infty$ as well. We require two more axioms:

(1) $\rho_n(x,y)<+\infty\Rightarrow\rho_{n+1}(x,y)=0\ (x,y\in X)$

(2) $\forall x,y\in X\ (x\ne y)\ \exists n\in\mathbb{Z}$ such that $0<\rho_n(x,y)<+\infty$

Then $\langle\{(x,y):x,y\in X,x\ne y\};\rho_n:n\in\mathbb{Z}\rangle$ is a small dimension structure.
\end{ex}

\P If $n<m,\rho_n(x,y)<+\infty$ then $\rho_m(x,y)=0$ by induction from (1).

Obviously $\langle\mathbb{Z},<\rangle$ is a complete ordered set that gives (ax2) and (ax3). 
\Pe

Here the dimension of an object (=2 points) gives the scale (or better to say name of the scale) where their distance can be reasonable measured. 

\bigskip

Now we construct an object that satisfies these axioms.

\begin{ex}Let 
$$X=\{(x_n)_{n=-\infty}^{+\infty}:x_n\in\mathbb{R}\text{ and }\exists k\in\mathbb{Z}\text{ such that }x_m=0\text{ whenever }m\geq k\}\subset\mathbb{R}^{\mathbb{Z}}.$$
For $x=(x_n),y=(y_n)\in X$ let $$\rho_n(x,y)=
\begin{cases}
|x_n-y_n|&\text{if }x_m=y_m\ \forall m>n\\
+\infty&\text{otherwise}
\end{cases}.$$

One can readily check that $\langle\{(x,y):x,y\in X,x\ne y\};\rho_n:n\in\mathbb{Z}\rangle$ satisfies the above axioms.
\end{ex}

We can define a continuous version as well. 
\begin{ex}\label{eca}Let $X$ be a set. For each $s\in\mathbb{R}$ let $\rho_s$ be a pseudo--metric on $X$ with the additional property that $\rho_s(x,y)$ can take the value $+\infty$ as well. We require two more axioms:

(1) If $x,y\in X,s,p\in S,s<p$ and $\rho_s(x,y)<+\infty$ then $\rho_{p}(x,y)=0$

(2) $\forall x,y\in X\ (x\ne y)\ \exists s\in\mathbb{R}$ such that $0<\rho_s(x,y)<+\infty$

Then $\langle\{(x,y):x,y\in X,x\ne y\};\rho_s:s\in\mathbb{R}\rangle$ is a dimension structure.
\end{ex}

\begin{ex}To get a slightly reduced such system let us take the dimension structure associated to the Hausdorff measure and dimension (Example \ref{eH}). For $0\leq s\leq 1,\ H,K\subset\mathbb{R}$ let $\rho_s(H,K)=\mu_s\big((H-K)\cup(K-H)\big)$. Evidently $\rho_s$ is a pseudo-metric and the system $\langle P(\mathbb{R}):\rho_s:s\in[0,1]\rangle$ satisfies condition (1)  in \ref{eca}. Of course $\rho_s$ is not defined for all $s\in\mathbb{R}$, just for $s\in[0,1]$ and (2) is not always holds.
\end{ex}

\begin{ex}Let $(X,\tau)$ be a topological space that is not $M_2$. Let a Borel measure $\mu$ be given on $X$ and let $\mu$ also denote the outer measure obtained by $\mu$. If $H\subset X$ let us define the accumulation points of $H$ regarding the system $\tau,\mu$:
$$H'=\{x\in H:\ K\text{ is open},x\in K\implies \mu(H\cap K)>0\}.$$
Set $H^{(0)}=H$ and $H^{(n+1)}=\big(H^{(n)}\big)'$ for $n\in\mathbb{N}$. Obviously $H^{(n+1)}\subset H^{(n)}$.

If $n\in\mathbb{N}\cup\{0\}$ let 
$$\mu_n(H)=
\begin{cases}
\mu\big(H^{(n)}\big)&\text{if }\mu\big(H^{(n+1)}\big)=0\\
+\infty&\text{otherwise}.
\end{cases}$$
Then $\langle P(X),\mu_n:n\in\mathbb{N}\cup\{0\}\rangle$  is a dimension structure.
\end{ex}

We can easily generalize some of the previous examples.

\begin{ex}Let $\langle X,<\rangle$ be a partially ordered set and $f:X\to X$ be a function with the property that if $x\in X$ then $f(x)\leq x$.  Let $f_0(x)=f(x)$ and $f_{n+1}(x)=f\big(f_n(x)\big)$ for $n\in\mathbb{N}$. If $n\in\mathbb{N}\cup\{0\}$ let 
$$\mu_n(x)=
\begin{cases}
0&\text{if }f_n(x)=f_{n+1}(x)\\
1&\text{if }f_n(x)\ne f_{n+1}(x)=f_{n+2}(x)\\
+\infty&\text{otherwise}.
\end{cases}$$
Then $\langle X,\mu_n:n\in\mathbb{N}\cup\{0\}\rangle$  is a dimension structure.
\end{ex}

\begin{ex}\label{eleb}If $H\subset\mathbb{R},H\in{\cal{L}}, n\in\mathbb{Z}$ let 
$$\mu_n(H)=\begin{cases}
+\infty&\text{if }\lambda(H\cap[n+1,+\infty))>0\\
\lambda(H\cap[n,n+1))&\text{otherwise}
\end{cases}$$
where $\lambda,{\cal{L}}$ denote the Lebesgue measure on $\mathbb{R}$ and the set of Lebesgue measurable sets respectively.

Then $\langle {\cal{L}},\mu_n: n\in\mathbb{Z}\rangle$ is a discrete dimension structure.
\end{ex}

We can create the planar version of the previous example.

\begin{ex}\label{eplleb}If $H\subset\mathbb{R}^2,H\in{\cal{L}}, n,m\in\mathbb{Z}$ let $$H_{n,m}=H\cap[n,n+1)\times[m,m+1),$$ 
$$\hat{H}_{n,m}=H\cap\{(x,y)\in\mathbb{R}^2:x\geq n+1\text{ or }y\geq m+1\},$$
and
$$\mu_{n,m}(H)=\begin{cases}
+\infty&\text{if }\lambda(\hat{H}_{n,m})>0\\
\lambda(H_{n,m})&\text{otherwise}
\end{cases}$$
where $\lambda$ denotes the Lebesgue measure on $\mathbb{R}^2$ and ${\cal{L}}$ is the set of Lebesgue measurabe sets.

Then $\big\langle {\cal{L}},\mu_{n,m}: (n,m)\in\langle\mathbb{Z}^2,<\rangle\big\rangle$ is a principal dimension structure where $<$ is the product order on $\mathbb{Z}^2$. Obviously here the underlying poset is not ordered.
\end{ex}
\P (ax1): Let $\mu_{n,m}(H)<+\infty,\ (n,m)<(n',m')$. Then $n\leq n',m<m'$ or $n<n',m\leq m'$ holds. Assume the first (the other is similar).

Evidently $\hat{H}_{n',m'}\subset \hat{H}_{n,m}$ which gives that $\lambda(\hat{H}_{n',m'})=0$ because $\lambda(\hat{H}_{n,m})=0$. Hence $\mu_{n',m'}(H)=\lambda(H_{n',m'})$ but $H_{n',m'}\subset \hat{H}_{n,m}$ which yields that $\mu_{n',m'}(H)=0$.

\smallskip

(ax2): Let $0<\mu_{n,m}(H)<+\infty,\ \mu_{n',m'}(H)<+\infty$. Then $\lambda(\hat{H}_{n,m})=0,\ \lambda(\hat{H}_{n',m'})=0$ and $0<\lambda(H_{n,m})<+\infty$. Therefore $H_{n,m}\not\subset \hat{H}_{n',m'}$. If $(n,m)$ and $(n',m')$ were incomparable then either $n<n',m>m'$ or $n>n',m<m'$ would hold. In any case $H_{n,m}\subset \hat{H}_{n',m'}$ would hold which is a contradiction.

\smallskip

(ax3): $\langle\mathbb{Z}^2,<\rangle$ is a discrete lattice therefore complete.

\smallskip

Principal: Let $(n,m)=\inf\{(k,l):\mu_{k,l}(H)<+\infty\}$ and $(n,m)<(n',m')$. First observe that $(n,m)\in \{(k,l):\mu_{k,l}(H)<+\infty\}$ because $\langle\mathbb{Z}^2,<\rangle$ is discrete. Hence $\mu_{n,m}(H)<+\infty$.

Suppose that $\mu_{n',m'}(H)=+\infty$. Then $\lambda(\hat{H}_{n',m'})>0$. But $\hat{H}_{n',m'}\subset \hat{H}_{n,m}$ which gives that $\lambda(\hat{H}_{n,m})>0$ and we would get that $\mu_{n,m}(H)=+\infty$ that is a contradiction.
\Pe

The following example is based on an expansion of numbers, namely every number $z$ in $[0,+\infty)$ can be uniquely written in the form 
$$z=e^{e^{\dots^{e^{\frac{x}{x+1}}}}}$$ where $x\in[0,+\infty)$. This form has two parameters: the height of the tower of $e^{e^{e^{\dots}}}$ (that can be 0) and $x$. Both parameters are uniquely determined by $z$.

\begin{ex}We need some definition first.

\medskip

Let $g_0:[0,+\infty)\to\mathbb{R},\ g_0(x)=\frac{x}{x+1}$. If $n\in\mathbb{N}$ then set $g_n:[0,+\infty)\to\mathbb{R},\ g_n(x)=e^{g_{n-1}(x)}$. Let us observe that if $n<m$ then $y\in Ran\ g_n,\ z\in Ran\ g_m$ implies that $y<z$. Moreover $\bigcup_{n\in\mathbb{N}\cup\{0\}}Ran\ g_n=[0,+\infty)$.

\medskip

If $x\in [0,+\infty),n\in\mathbb{N}\cup\{0\}$ then set $h_{x,n}(y):[0,+\infty]\to\mathbb{R}$
$$h_{x,n}(y)=
\begin{cases}
g_n(y)-x&\text{if }y<+\infty\\
\lim\limits_{z\to+\infty}g_n(z)-x&\text{if }y=+\infty.
\end{cases}$$
The limit exists because $z\mapsto g_n(z)-x$ strictly increasing. Hence $h_{x,n}(y)$ strictly increasing.

\medskip

Let us show that if $x\in[0,+\infty),n\in\mathbb{N}\cup\{0\}$ are given then there is a unique $y_0\in[0,+\infty]$ such that $|h_{x,n}(y)|$ takes its minimum in $y_0$.

The minimum exists because $|h_{x,n}(y)|$ is continuous and $[0,+\infty]$ is compact. The minimum is unique since $h_{x,n}(y)$ strictly increasing.

\medskip

Now we can define the dimension structure.

Let ${\cal{D}}=\langle [0,+\infty),\mu_s:s\in \mathbb{N}\cup\{0\}\rangle$ where both $[0,+\infty)$ and $\mathbb{N}\cup\{0\}$ is equipped with the usual order and if $x\in[0,+\infty),n\in\mathbb{N}\cup\{0\}$ then $\mu_n(x)=y_0$ where $y_0$ is the unique number in $[0,+\infty]$ such that $|h_{x,n}(y)|$ takes its minimum.
Then ${\cal{D}}$ is a finitely synchronized normal dimension structure.
\end{ex}

\P We have to prove (ax1'). Suppose that $\mu_n(x)<+\infty$ for some $x\in[0,+\infty),n\in\mathbb{N}\cup\{0\}$. It is easy to see that then either $x\in Ran\ g_n$ or $\forall z\in Ran\ g_n,\ x<z$ holds. Both cases imply that $\forall z\in Ran\ g_{n+1},\ x<z$ holds which yields that $\mu_{n+1}(x)=0$.

Finite synchronization follows from \ref{psynf} and \ref{psync1} and the fact that $x\leq y$ implies that $\forall s\in S\ \mu_s(x)\leq\mu_s(y)$. 

Normality is trivial.
\Pe

Here the dimension of a point $x$ is the height of the tower of $e^{e^{e^{\dots}}}$.

\smallskip

A straightforward examle can show that ${\cal{D}}$ is not $\aleph_0$ synchronized.

\section{Building new dimension structures}

There are many ways how we can get a new dimension structure from already existing ones. In this section we are going to present such constructions.

Most of the proofs are straightforward therefore they are omitted in many cases.

\subsection{Substructure}

\begin{prp}Let ${\cal{D}}=\langle X,\mu_s:s\in S\rangle$ be a dimension structure. Let $Y\subset X, P\subset S$ such that $P$ inherits the order from $S$ and $P$ is complete. When $p\in P$ then let $\nu_p=\mu_p|_Y$. 

Then ${\cal{D}}'=\langle Y,\nu_p:p\in P\rangle$ is a dimension structure. \Pes
\end{prp}

\begin{df}Let ${\cal{D}}=\langle X,\mu_s:s\in S\rangle,\ {\cal{D}}'=\langle Y,\nu_p:p\in P\rangle$ be two dimension structures. Let $Y\subset X, P\subset S$ such that $P$ inherits the order from $S$. When $p\in P$ then let $\nu_p=\mu_p|_Y$. 

Then ${\cal{D}}'$ is called the substructure of ${\cal{D}}$.
\end{df}

This method contains the important case when we leave out measures only and keep $X$ intact. Even more important case is defined in
\begin{df}Let ${\cal{D}}=\langle X,\mu_s:s\in S\rangle$ be a dimension structure. Set $X'=\{x\in X:\exists s\in S$ such that $0<\mu_s(x)<+\infty\}$ and $S'=\{s\in S:\exists x\in X$ such that $0<\mu_s(x)<+\infty\}$.
${\cal{D}}'=\langle X',\mu_s:s\in S'\rangle$ is called the normalization of ${\cal{D}}$ and it is a normal dimension structure.
\end{df}

\begin{rem}One can readily create examples where the dimension of a point gets smaller in the substructure, and also where the dimension of a point gets greater in the substructure than it was in the original structure.

E.g. let $\langle P(\mathbb{R}):\mu_s:s\in[0,1]\rangle$ as in example \ref{eH}. Let $Y=X=P(\mathbb{R})$ and $P=[0,0.3]\cup(0.5,0.6]\cup[0.8,1]$. Let ${\cal{D}}'=\langle Y,\nu_p:p\in P\rangle$.
Let $H,K\subset\mathbb{R}$ such that the Hausdorff dimension of $H,K$ equals to 0.4 and 0.7 respectively. Then $\dim_{{\cal{D}}'} H=0.3$ and $\dim_{{\cal{D}}'} K=0.8$.
\end{rem}

\begin{prp}Let ${\cal{D}}'$ be a substructure of ${\cal{D}}$ as above. Let $x\in Y$ be an $s$-point in ${\cal{D}}$ for $s\in P$. Then $x$ is an $s$-point in ${\cal{D}}'$ as well. I.e. $\dim_{{\cal{D}}'} x=\dim_{{\cal{D}}} x$ in this case.
\end{prp}

\subsection{Quotient structure}

\begin{prp}Let ${\cal{D}}=\langle X,\mu_s:s\in S\rangle$ be a dimension structure with $S$ being complete. Let $X=\cup_{i\in I}X_i, X_i\cap X_j=\emptyset\ (i,j,\in I,i\ne j)$. For $i\in I,s\in S$ let $\nu_s(i)=\sup\{\mu_s(x):x\in X_i\}$

Then ${\cal{D}}'=\langle I,\nu_s:s\in S\rangle$ is a dimension structure.
\end{prp}

\P Let us check (ax1). If $\nu_s(i)<+\infty,s<p\in S$ then $\mu_s(x)<+\infty\ \forall x\in X_i$. It gives that $\mu_p(x)=0\ \forall x\in X_i$ that implies that $\nu_p(i)=0$.

Now we show (ax2).  If $0<\nu_s(i)<+\infty$ then $\exists x_0\in X_i$ such that $0<\mu_s(x_0)<+\infty$. If $\nu_p(i)<+\infty$ then $\mu_p(x)<+\infty\ \forall x\in X_i$. Hence $\mu_p(x_0)<+\infty$ and by (ax2) in ${\cal{D}}$ we get that $s,p$ are comparable.
\Pes

\begin{df}${\cal{D}}'$ is called the quotient of ${\cal{D}}$.
\end{df}

\begin{ex}(ax3) does not necessarily hold if $S$ is not complete.

Let $S=\{-\frac{1}{n}:n\in\mathbb{N}\}\cup\{\frac{1}{n}:n\in\mathbb{N}\}$ with the order inherited from $\mathbb{R}$. Let $X=\{x_1,x_2,x_3,\dots\}$ and if $s\in S, x_i\in X$ then let 
$$\mu_s(x_i)=
\begin{cases}
0&\text{if }s\geq -\frac{1}{i}\\
+\infty&\text{otherwise.}
\end{cases}$$
Now let $X=X_1$ a one element partition (i.e. $I=\{1\}$). We get that
$$\nu_s(1)=
\begin{cases}
0&\text{if }s>0\\
+\infty&\text{otherwise.}
\end{cases}$$
Therefore there is no infimum in $S$ of the set $\{s\in S:\nu_s(1)=0\}=\{\frac{1}{n}:n\in\mathbb{N}\}$. 
\end{ex}

We give a generic example for a quotient.

\begin{prp}Let ${\cal{D}}=\langle X,\mu_s:s\in S\rangle$ be a dimension structure with $S$ being complete. Let us take the partition of $X$
$$X=\bigcup\limits_{(d,m)\in \bar{S}\times[0,+\infty]}C^{\cal{D}}_{d,m}$$
and set $Y=\{(d,m)\in \bar{S}\times[0,+\infty]:\ \exists x\in X\text{ such that }\dim x=d,\mu_d(x)=m\}$ and $\nu_d(y)=m\ (y\in Y,y=(d,m))$. Then $\langle Y,\nu_s:s\in S\rangle$ is a quotient of ${\cal{D}}$. \Pes
\end{prp}

\begin{prp}Let ${\cal{D}},{\cal{D}}'$ be dimension structures as above with $S$ being complete. Then $\sup\{\dim_{{\cal{D}}} x:x\in X_i\}\leq\dim_{{\cal{D}}'} i$.
\end{prp}
\P Let $s=\dim_{{\cal{D}}'} i$ and $\nu_s(i)=m$. Let us suppose that there is $y\in X_i$ such that $s<\dim_{{\cal{D}}} y=k$. Then $\mu_s(y)=+\infty$ which gives that $m=+\infty$ too. Also $\nu_k(i)=0$ which gives that $\mu_k(y)=0$. From $s=\dim_{{\cal{D}}'} i<k$ we get that there is $p\in S$ such that $s<p<k$. Then $\nu_p(i)=0$ and $\mu_p(y)=+\infty$ holds which yields that $\nu_p(i)=+\infty$ which is a contradiction.
\Pes

\subsection{Sum of structures}

\begin{prp}Let $P$ be a partially ordered set and for all $p\in P$ let ${\cal{D}}_p=\langle X_p,\mu^{(p)}_s:s\in S_p\rangle$ be a dimension structures such that $X_p\cap X_q=\emptyset$ and $S_p\cap S_q=\emptyset\ (p,q\in P,p\ne q)$. Let $\pm\infty_{S_p}\in S_p$ denote the minimum and maximum of $S_p$.

Let $X=\bigcup\limits_{p\in P}X_p, S=\bigcup\limits_{p\in P}S_p$. Let us define partial order on $S$: if $s_1,s_2\in S,\ s_1\in S_{p_1},s_2\in S_{p_2}\ (p_1,p_2\in P)$ then set 
$$s_1<s_2\iff\begin{cases}
\text{if }p_1=p_2\text{ and }s_1<s_2\text{ in }{S_{p_1}}\\
\text{if }p_1<p_2.
\end{cases}$$
For $x\in X, s\in S,\ x\in X_p,s\in S_q\ (p,q\in P)$ let 
$$\eta_s(x)=\begin{cases}
\mu_s(x)&\text{if }p=q\\
0&\text{if }p<q\\
+\infty&\text{if }q<p\\
+\infty&\text{if }p,q\text{ are non-comparable}.
\end{cases}$$

Then ${\cal{D}}'=\langle X,\eta_s:s\in S\rangle$ is a dimension structure. \Pes
\end{prp}

\begin{df}${\cal{D}}'$ is called the sum of $\{{\cal{D}}_p:p\in P\}$ over $P$ and denoted by $\oplus_{p\in P}{\cal{D}}_p$.
\end{df}

There is one slight disadvantage of this method: if for a point in $X_p$ the dimension was $+\infty_{S_p}$ in the original structure then it may significantly change in the new structure.

\begin{prp}Let $P$ be a partially ordered set and for all $p\in P$ let ${\cal{D}}_p=\langle X_p,\mu^{(p)}_s:s\in S_p\rangle$ be a dimension structures as above. Let ${\cal{D}}=\oplus_{p\in P}{\cal{D}}_p$
If $x\in X, x\in X_p$ then
$$\dim_{\cal{D}} x=
\begin{cases}
\dim_{{\cal{D}}_p} x&\text{if }\dim_{{\cal{D}}_p} x\ne+\infty_{S_p}\\
\dim_{{\cal{D}}_p} x&\text{if }\dim_{{\cal{D}}_p} x=+\infty_{S_p}\text{ and }\mu_{+\infty_{p}}(x)<+\infty\\
\dim_{{\cal{D}}_p} x&\text{if }\dim_{{\cal{D}}_p} x=+\infty_{S_p}\text{ and }\mu_{+\infty_{p}}(x)=+\infty\text{ and }\not\exists\min\{q\in P:p<q\}\\
-\infty_{S_r}&\text{if }\dim_{{\cal{D}}_p} x=+\infty_{S_p}\text{ and }\mu_{+\infty_{p}}(x)=+\infty\text{ and }r=\min\{q\in P:p<q\}.
\end{cases}$$
\Pes
\end{prp}

\subsection{Measure sum}

\begin{prp}Let ${\cal{D}}_i=\langle X,\mu^{(i)}_s:s\in S\rangle$ be dimension structures for all $i\in\mathbb{N}$ with $S$ being complete. Then $${\cal{D}}=\Big\langle X,\sum\limits_{i=1}^{\infty}\mu^{(i)}_s:s\in S\Big\rangle$$ is a dimension structure. 
\end{prp}
\P (ax1): Let $\mu_s=\sum\limits_{i=1}^{\infty}\mu^{(i)}_s$ and $x\in X,s,p\in S,s<p$.
If $\mu_s(x)<+\infty$ then $\forall i\in\mathbb{N}\ \mu^{(i)}_s(x)<+\infty$ which implies that $\forall i\ \mu^{(i)}_p(x)=0$ that gives that $\mu_p(x)=0$.

(ax2): If $0<\sum\limits_{i=1}^{\infty}\mu^{(i)}_s(x)<+\infty$ then there is $j\in\mathbb{N}$ such that $0<\mu^{(j)}_s(x)<+\infty$. If $\sum\limits_{i=1}^{\infty}\mu^{(i)}_p(x)<+\infty$ then $\forall i\in\mathbb{N}\ \mu^{(i)}_p(x)<+\infty$ holds, especially $\mu^{(j)}_p(x)<+\infty$. That gives that $s,p$ are comparable.
\Pes

\begin{df}We call ${\cal{D}}$ the measure sum of the system $\{{\cal{D}}_i:i\in\mathbb{N}\}$.
\end{df}

\begin{ex}\label{exmsnc}If $S$ is not complete then (ax3) does not necessarily hold.

Let $X=\{x\},\ S'=\{\pm\frac{1}{n}:n\in\mathbb{N}\},\ S=S'\cup\{a\}$. If $s,p\in S$ then let 
$$s<p\iff
\begin{cases}
s,p\in S'\text{ and }s<p\\
s=a,p=1\\
s=-1,p=a.
\end{cases}$$
Evidently $S$ is a lattice. Let $$\mu^{(1)}_s(x)=
\begin{cases}
0&\text{if }s\in S',s\not=-1\\
+\infty&\text{if }s=-1\text{ or }s=a,
\end{cases}$$
and
$$\mu^{(2)}_s(x)=
\begin{cases}
0&\text{if }s\in S',s>0\text{ or }s=a\\
+\infty&\text{otherwise}.
\end{cases}$$
Then obviously there does not exists the infimum of $\{s\in S:\mu^{(1)}_s(x)+\mu^{(2)}_s(x)=0\}=\{\frac{1}{n}:n\in\mathbb{N}\}$.
\end{ex}

Obviously we can have the measure sum of finitely many dimension structures too. For principal dimension structures if $S$ is a lattice then for finitely many dimension structures we can omit the condition that $S$ is complete.

\begin{prp}Let ${\cal{D}}_i=\langle X,\mu^{(i)}_s:s\in S\rangle$ be principal dimension structures for $1\leq i\leq n\ (n\in\mathbb{N})$ with $S$ being a lattice. Then $${\cal{D}}=\Big\langle X,\sum\limits_{i=1}^{n}\mu^{(i)}_s:s\in S\Big\rangle$$ is a dimension structure.
\end{prp}
\P We have to check (ax3).  Let $S'=\{s\in S:\sum\limits_{i=1}^{n}\mu^{(i)}_s(x)<+\infty\}, S'_i=\{s\in S:\mu^{(i)}_s(x)<+\infty\}\ (1\leq i\leq n)$. Clearly $S'=\bigcap\limits_{i=1}^{n}S'_i$.
By principality condition $S'_i\cup\{s'_i\}$ is a principal filter ($1\leq i\leq n$) where $s'_i=\inf S'_i$. In a lattice the intersection of finitely many principal filters is a principal filter moreover $s'=\sup\{s'_i:1\leq i\leq n\}=\inf \bigcap_{i=1}^n (S'_i\cup\{s'_i\})$. We show that $\inf S'=s'$ too. There are two cases.
\begin{itemize}
\item[(1)] If $\exists i\in\{1,\dots,n\}$ such that $s'=s'_i$. Then clearly $S'=S'_i$ and the statement is obvious.
\item[(2)] If no such $i$ exists then $s'>s'_i$ for all $i\in\{1,\dots,n\}$. That means that $s'\in S'_i$ for all $i$, i.e. $s'\in\bigcap_{i=1}^n S'_i=S'$ and evidently it is the smallest element in $S'$.
\end{itemize}
\Pes

\begin{ex}We cannot ommit the principality condition as Example \ref{exmsnc} shows.

\end{ex}

\begin{prp}\label{pmsdcom}Let ${\cal{D}}$ be the measure sum of the dimension structures ${\cal{D}}_i=\langle X,\mu^{(i)}_s:s\in S\rangle\ (i\in\mathbb{N})$ ($S$ is complete). If $x\in X$ then $\sup\{\dim_{{\cal{D}}_i} x:i\in\mathbb{N}\}\leq\dim_{{\cal{D}}} x$. Moreover if $S$ is dense then $\sup\{\dim_{{\cal{D}}_i} x:i\in\mathbb{N}\}=\dim_{{\cal{D}}} x$.
\end{prp}
\P Set $d=\dim_{{\cal{D}}} x,e=\sup\{\dim_{{\cal{D}}_i} x:i\in\mathbb{N}\}$. 

If $d< s\in S$ then $\mu_s(x)=0$ which implies that $\mu^{(i)}_s(x)=0\ (\forall i\in\mathbb{N})$ hence $\dim_{{\cal{D}}_i} x\leq d$ i.e. $e\leq d$.

Let us now assume that $e<d$ for some $x\in X$ while $S$ being dense. Then there is $p\in S$ such that $e<p<d$. We get that $\mu^{(i)}_p(x)=0\ (i\in\mathbb{N})$ which yields that $\mu_p(x)=0$ and then $d\leq p$ which is a contradiction.
\Pes

\begin{rem}In the previous proposition we cannot omit the density for $S$.

For showing that let ${\cal{D}}_i$ equal to the same dimension structure for all $i\in\mathbb{N}$ with $S$ being discrete. Let $x\in X,s\in S$ such that $0<\mu^{(i)}_s(x)<+\infty$. Then $\mu^{(i)}_{s^+}(x)=0\ (i\in\mathbb{N})$ hence $\mu_{s^+}(x)=0$ and $\mu_{s}(x)=+\infty$. Therefore $\sup\{\dim_{{\cal{D}}_i} x:i\in\mathbb{N}\}=s$ while $\dim_{{\cal{D}}} x=s^+$.
\end{rem}

\begin{prp}Let ${\cal{D}}$ be the measure sum of the dimension structures ${\cal{D}}_i=\langle X,\mu^{(i)}_s:s\in S\rangle\ (i\in\mathbb{N})$ with $S$ being discrete ordered. If $x\in X$ let $d_x=\sup\{\dim_{{\cal{D}}_i} x:i\in\mathbb{N}\}$. If $d_x<\dim_{{\cal{D}}} x$ then $(d_x)^+=\dim_{{\cal{D}}} x$.
\end{prp}
\P Obviously $\mu^{(i)}_{(d_x)^+}(x)=0\ (i\in\mathbb{N})$ which implies that $\mu_{(d_x)^+}(x)=0$ i.e. $\dim_{{\cal{D}}} x\leq (d_x)^+$ and \ref{pmsdcom} gives the statement.
\Pe

We show another such construction.

\begin{prp}Let ${\cal{D}}_i=\langle X,\mu^{(i)}_s:s\in S\rangle$ be dimensional structures for all $i\in I$ with $S$ being complete. Then ${\cal{D}}=\big\langle X,\sup\{\mu^{(i)}_s:i\in I\}:s\in S\big\rangle$ is a dimensional structure. \Pes
\end{prp}

\subsection{Direct product}

\begin{prp}Let ${\cal{D}}_1=\langle X_1,\mu^{(1)}_s:s\in S_1\rangle,\ {\cal{D}}_2=\langle X_2,\mu^{(2)}_s:s\in S_2\rangle$ be two dimension structures.
Let $X=X_1\times X_2,S=S_1\times S_2$ equipped with the product order.
For $x\in X,x=(x_1,x_2),s\in S,s=(s_1,s_2)$ let 
$$\eta_s(x)=
\begin{cases}
+\infty&\text{if }\mu^{(1)}_{s_1}(x_1)=+\infty\text{ or }\mu^{(2)}_{s_2}(x_2)=+\infty\\
\mu^{(1)}_{s_1}(x_1)\cdot\mu^{(2)}_{s_2}(x_2)&\text{otherwise}.
\end{cases}$$

Then ${\cal{D}}'=\langle X,\eta_s:s\in S\rangle$ is a dimension structure.
\end{prp}

\P (ax1): Let $s,s'\in S,x\in X,s=(s_1,s_2),x=(x_1,x_2),s<s'=(s_1',s_2')$. Then either $s_1<s'_1,s_2\leq s_2'$ or $s_2<s_2',s_1\leq s_1'$. 

If $\eta_s(x)<+\infty$ then $\mu^{(1)}_{s_1}(x_1)<+\infty$ and $\mu^{(2)}_{s_2}(x_2)<+\infty$. There are two cases.

1. If $s_1<s_1'$ then $\mu^{(1)}_{s_1'}(x)=0$ hence $\eta_{s'}(x)=0$.

2. If $s_2<s_2'$ then $\mu^{(2)}_{s_2'}(x_2)=0$ which implies again that $\eta_{s'}(x)=0$.

\smallskip

(ax2): Let $0<\eta_s(x)<+\infty,\eta_{s'}(x)<+\infty$ and $x=(x_1,x_2)$ and $s=(s_1,s_2),s'=(s_1',s_2')$. 
Then we know that $0<\mu^{(1)}_{s_1}(x_1)<+\infty$ and $0<\mu^{(2)}_{s_2}(x_2)<+\infty$, $\mu^{(1)}_{s_1'}(x_1)<+\infty$ and $\mu^{(2)}_{s_2'}(x_2)<+\infty$.
Then $s_1\leq s_1'$ and $s_2\leq s_2'$ which gives that $s\leq s'$.

\smallskip

(ax3): If $x\in X,x=(x_1,x_2)$ then set $S'=\{s\in S:\eta_s(x)<+\infty\},S_1'=\{s_1\in S_1:\mu^{(1)}_{s_1}(x_1)<+\infty\},S_2'=\{s_2\in S_2:\mu^{(2)}_{s_2}(x_2)<+\infty\}$. Then $S'=S_1'\times S_2'$ and $\inf S'=(\inf S_1',\inf S_2')$ hence $\inf S'$ exists.
\Pes

\begin{df}We call ${\cal{D}}'$ the direct product of ${\cal{D}}_1$ and ${\cal{D}}_2$ and denote it by ${\cal{D}}_1\times{\cal{D}}_2$.
\end{df}

\begin{rem}Actually we showed that $\dim_{{\cal{D}}_1\times{\cal{D}}_2}(x_1,x_2)=(\dim_{{\cal{D}}_1} x_1,\dim_{{\cal{D}}_2} x_2)$.
\end{rem}

\begin{rem}In Examples \ref{eleb} and \ref{eplleb} let us restrict ourself to sets with positive measure. 
Then the direct product of the structure in Example \ref{eleb} with itself can be embedded to a substructure of the structure in Example \ref{eplleb} if we identify $(H,K)$ with $H\times K\ (H,K\subset\mathbb{R})$.
\end{rem}

\begin{rem}Let $f:(\mathbb{R}^+\cup\{0\})\times(\mathbb{R}^+\cup\{0\})\to\mathbb{R}^+$ and $f(0,x)=f(x,0)=0\ (\forall x\in\mathbb{R}^+\cup\{0\})$.
If in the definition of $\eta_s(x)$ we replace the product (of two numbers) to $f$ then we end up with a dimension structure as well. 
\end{rem}

\begin{prp}The direct product of two p-small dimension structures is p-small. \Pes
\end{prp}

\begin{prp}Let ${\cal{D}}_1=\langle X_1,\mu^{(1)}_s:s\in S_1\rangle,\ {\cal{D}}_2=\langle X_2,\mu^{(2)}_s:s\in S_2\rangle$ be two p-small, principal dimension structures.
Then ${\cal{D}}_1\times{\cal{D}}_2$ is principal as well.
\end{prp}
\P Let $x=(x_1,x_2)\in X_1\times X_2$, $S'=\{(s_1,s_2)\in S_1\times S_2:\mu_{(s_1,s_2)}(x_1,x_2)<+\infty\}$, $S'_i=\{s_i\in S_i:\mu^{(i)}_{s_i}(x_i)<+\infty\}\ (i=1,2)$. Let $s'_i=\inf S'_i\ (i=1,2)$. Then we know that $\inf S'=(s'_1,s'_2)$. By assumption $\mu_{s'_1(x_1)}<+\infty,\ \mu_{s'_2(x_2)}<+\infty$. 

If $(s'_1,s'_2)<(s''_1,s''_2)\in S_1\times S_2$ then $\exists i\in\{1,2\}$ such that $s'_i<s''_i$ that gives $\mu_{s''_i}(x_i)=0$ which implies that $\mu_{(s''_1,s''_2)}(x_1,x_2)=0$ that we had to prove by \ref{pprinceq}.
\Pe

Now we define the i-direct product of arbitrary many structures.

\begin{prp}Let ${\cal{D}}_i=\langle X_i,\mu_s:s\in S_i\rangle$ be dimension structures for $i\in I$.
Let $Z=\bigtimes\limits_{i\in I} X_i,\ Q=\bigtimes\limits_{i\in I} S_i$ equipped with the product order.
Let $z\in Z,q\in Q$ and let us denote the $i^{th}$ coordinate by $z_i,q_i$ respectively. Let 
$$\eta_q(z)=
\begin{cases}
+\infty&\text{if }\exists i\in I\text{ such that }\mu_{q_i}(z_i)=+\infty\\
\inf\{\mu_{q_i}(z_i):i\in I\}&\text{otherwise}.
\end{cases}$$

Then ${\cal{D}}'=\langle Z,\eta_q:q\in Q\rangle$ is a dimension structure.
\end{prp}

\P (ax1): Let $q,q'\in Q,z\in Z,q<q'$. Then $q_j\leq q'_j$ and there is $i\in I$ such that $q_i< q'_i$. 

If $\eta_q(z)<+\infty$ then $\mu_{q_j}(z_j)<+\infty\ \forall j\in I$. Hence $\forall j\in I\ \mu_{q'_j}(z_j)<+\infty$ moreover $\mu_{q'_i}(z_i)=0$. Therefore $\eta_{q'}(z)=0$.

\smallskip

(ax2): Let $0<\eta_q(z)<+\infty,\eta_{q'}(z)<+\infty$. 
Then we know that $0<\mu_{q_j}(z_j)<+\infty$ and $\mu_{q'_j}(z_j)<+\infty$ $(\forall j\in I)$. Which implies that $\forall j\in I\ q_j\leq q'_j$ which gives that $q\leq q'$.

\smallskip

(ax3): If $z\in Z$ then set $Q'=\{q\in Q:\eta_q(z)<+\infty\},S'_i=\{s\in S_i:\mu_s(z_i)<+\infty\}$ $(i\in I)$. Then $Q'=\bigtimes\limits_{i\in I} S'_i$ hence $\inf Q'$ exists.
\Pes

\begin{df}We call ${\cal{D}}'$ the i-direct product of the system $({\cal{D}}_j:j\in I)$ and denote it by ${i \atop \bigtimes}_{j\in I}{\cal{D}}_j$.
\end{df}

\begin{rem}For i-direct product $(\dim_{{\cal{D}}'}z)_i=\dim_{{\cal{D}}_i}z_i$ holds.
\end{rem}

Now we define another type of product.

\begin{prp}\label{plpm}Let ${\cal{D}}_1=\langle X_1,\mu^{(1)}_s:s\in S_1\rangle,\ {\cal{D}}_2=\langle X_2,\mu^{(2)}_s:s\in S_2\rangle$ be two dimension structures and let ${\cal{D}}_1$ be small.

Let $X=X_1\times X_2,S=S_1\times S_2$ equipped with the lexicographic order.
For $x\in X,x=(x_1,x_2),s\in S,s=(s_1,s_2)$ let 
$$\mu_{s}(x)=\mu_{(s_1,s_2)}(x_1,x_2)=
\begin{cases}
+\infty&\text{if }\mu_{s_1}(x_1)=+\infty\\
\mu_{s_1}(x_1)\cdot\mu_{s_2}(x_2)&\text{otherwise}.
\end{cases}$$

Then ${\cal{D}}'=\langle X,\mu_s:s\in S\rangle$ is a dimension structure.
\end{prp}
\P \begin{itemize}
\item[(ax1):]If $\mu_{(s_1,s_2)}(x_1,x_2)<+\infty$ then either $\mu_{s_1}(x_1)=0$ or $0<\mu_{s_1}(x_1)<+\infty$ and $\mu_{s_2}(x_2)<+\infty$.

Let $(s_1,s_2)<(s'_1,s'_2)$. There are two cases.

(1) $s_1=s'_1,s_2<s'_2$: If $\mu_{s_1}(x_1)=0$ then $\mu_{s'_1}(x_1)=0$ hence $\mu_{s}(x)=0$. If $\mu_{s_2}(x_2)<+\infty$ then $\mu_{s'_2}(x_2)=0$ therefore $\mu_{s}(x)=0$.

(2) $s_1<s'_1$: Then $\mu_{s'_1}(x_1)=0$ hence $\mu_{s}(x)=0$.
\item[(ax2):]If $0<\mu_{(s_1,s_2)}(x_1,x_2)<+\infty$ then $0<\mu_{s_1}(x_1)<+\infty$ and $0<\mu_{s_2}(x_2)<+\infty$. Let $\mu_{(s'_1,s'_2)}(x_1,x_2)<+\infty$. Then either $\mu_{s'_1}(x_1)=0$ or $0<\mu_{s'_1}(x_1)<+\infty$ and $\mu_{s'_2}(x_2)<+\infty$. When applying (ax2) in ${\cal{D}}_1$ and ${\cal{D}}_2$ we get that $s_1<s_1'$ or $s_1= s'_1,s_2\leq s_2'$. That gives that $(s_1,s_2)\leq(s'_1,s'_2)$.
\item[(ax3):]For a given $x=(x_1,x_2)\in X$ let $$S'=\{s\in S:\mu_s(x)=0\},\ S'_i=\{s_i\in S_i:\mu_{s_i}(x_i)=0\}\ (i=1,2).$$
Then using the fact that $x_1$ is a dim-point we get that $$S'=S'_1\times S_2\cup\{\dim_{{\cal{D}}_1} x_1\}\times S'_2$$
which implies that $\dim(x_1,x_2)=(\dim x_1,\dim x_2)$.
\end{itemize}
\Pes

\begin{df}We call ${\cal{D}}'$ the l-direct product of ${\cal{D}}_1$ and ${\cal{D}}_2$ and denote it by ${\cal{D}}_1{l \atop \bigtimes}{\cal{D}}_2$.
\end{df}

\begin{prp}Let ${\cal{D}}_1=\langle X_1,\mu^{(1)}_s:s\in S_1\rangle,\ {\cal{D}}_2=\langle X_2,\mu^{(2)}_s:s\in S_2\rangle$ be two principal dimension structures and let ${\cal{D}}_1$ be p-small.
Then ${\cal{D}}_1{l \atop \bigtimes}{\cal{D}}_2$ is principal as well.
\end{prp}
\P For a given $x=(x_1,x_2)\in X$ let $$S'=\{s\in S:\mu_s(x)<+\infty\},\ S'_i=\{s_i\in S_i:\mu_{s_i}(x_i)<+\infty\}\ (i=1,2).$$
As in the proof of the previous proposition (\ref{plpm}) we get that $$S'=S'_1\times S_2\cup\{\dim_{{\cal{D}}_1} x_1\}\times S'_2.$$
Set $s_i=\inf S'_i\ (i=1,2)$. With that we have that $\dim x=(s_1,s_2)$.

Let $s'=(s'_1,s'_2)\in S$ such that $(s_1,s_2)<(s'_1,s'_2)$. There are two cases.
\begin{itemize}
\item[(1)]If $s_1<s'_1$ then $\mu_{s'_1}(x_1)=0$ which yields that $\mu_{s'}(x)=0$.
\item[(2)]If $s_1=s'_1,s_2<s'_2$ then by p-smallness $\mu_{s'_1}(x_1)<+\infty$ and $\mu_{s'_2}(x_2)=0$ which gives again that $\mu_{s'}(x)=0$.
\end{itemize}
By \ref{pprinceq} we are done.
\Pes

\subsection{Mapping between structures}

\begin{prp}Let ${\cal{D}}_1=\langle X_1,\mu^{(1)}_s:s\in S\rangle$ be a dimension structure with $S$ being complete ordered. Let ${\cal{D}}_2=\langle X_2,\mu^{(2)}_s:s\in S\rangle$ be given such that $X_2$ is a set, $\mu^{(2)}_s\ (s\in S)$ is a function $\mu^{(2)}_s:X_2\to\mathbb{R}^+\cup\{0,+\infty\}$. 
Let $f:X_1\to X_2$ be a surjective mapping. Let $\forall x\in X_1\ \mu_s^{(1)}(x)\leq\mu^{(2)}_s(f(x))$ and if $\forall x\in f^{-1}(y)\ (y\in X_2)\ \mu_s^{(1)}(x)=0$ then $\mu^{(2)}_s(y)=0$. Then ${\cal{D}}_2$ is a dimension structure as well.
\end{prp}

\P If $\mu^{(2)}_s(y)<+\infty$ then there is $x\in X$ such that $y=f(x)$ and then $\mu_s^{(1)}(x)<+\infty$ holds as well. If $s<p$ then $\mu_p^{(1)}(x)=0$ which gives that $\mu^{(2)}_p(y)=0$ too.
\Pe

The previous condition holds if $\sup\{\mu_s^{(1)}(x):f(x)=y\}=\mu^{(2)}_s(y)$.

\begin{rem}If we substitute the condition 
$$\forall x\in X_1\ \mu_s^{(1)}(x)\leq\mu^{(2)}_s(f(x))$$  to 
$$\forall x\in X_1\ sign(\mu_s^{(1)}(x))=sign(\mu^{(2)}_s(f(x)))$$ we get a valid statemet too.
\end{rem}

\begin{df}Let ${\cal{D}}_i=\langle X_i,\mu^{(i)}_s:s\in S_i\rangle$ be two dimension structures for $i=1,2$. Let $f:X_1\to X_2,\ g:S_1\to S_2$ be given such that $g$ preserves the order (if $s\leq p\ (s,p\in S_1)$ then $g(s)\leq g(p)$). 

We call $(f,g)$ a morphism between ${\cal{D}}_1$ and ${\cal{D}}_2$ if $x\in X_1,s\in S_1$ implies that $\mu^{(1)}_s(x)\leq \mu^{(2)}_{g(s)}(f(x))$.

We call $(f,g)$ an isomorphism between ${\cal{D}}_1$ and ${\cal{D}}_2$ if $f,g$ are bijections and $x\in X_1,s\in S_1$ implies that $\mu^{(1)}_s(x)= \mu^{(2)}_{g(s)}(f(x))$.

We call $(f,g)$ a semi-isomorphism between ${\cal{D}}_1$ and ${\cal{D}}_2$ if $f,g$ are bijections and $x\in X_1,s\in S_1$ implies that $$\mu^{(1)}_s(x)<+\infty\iff \mu^{(2)}_{g(s)}(f(x))<+\infty.$$
\end{df}

\begin{prp}If $(f,g)$ is a semi-isomorphism between ${\cal{D}}_1$ and ${\cal{D}}_2$ then $$\mu^{(1)}_s(x)=+\infty\iff \mu^{(2)}_{g(s)}(f(x))=+\infty\ (x\in X,s\in S).$$
\end{prp}

\begin{prp}If $(f,g)$ is a semi-isomorphism between ${\cal{D}}_1$ and ${\cal{D}}_2$ then  $g(\dim x)=\dim f(x)\ (x\in X)$.

If $(f,g)$ is a isomorphism between ${\cal{D}}_1$ and ${\cal{D}}_2$ then $$\mu^{(1)}_{\dim x}(x)=\mu^{(2)}_{g(\dim x)}(f(x))\ (x\in X)$$ holds as well. \Pes
\end{prp}

\begin{lem}If $(f,g)$ is a morphism between ${\cal{D}}_1$ and ${\cal{D}}_2$ and $0<\mu^{(1)}_{\dim x}(x)$  and $g(\dim x)$ and $\dim f(x)$ are comparable
then  $g(\dim x)\leq\dim f(x)\ (x\in X)$.
\end{lem}
\P From the condition we get that $0<\mu^{(2)}_{g(\dim x)}(f(x))$ which yields the statement.
\Pes

\begin{prp}Let ${\cal{D}}_i=\langle X_i,\mu^{(i)}_s:s\in S_i\rangle$ be two dimension structures ($i=1,2$) with $S_1,S_2$ being dense, complete and ordered. Let $f:X_1\to X_2,\ g:S_1\to S_2$ be given such that $g$ preserves the order and let $g$ be continuous between the topologies induced by the orders. Then $g(\dim x)\leq\dim f(x)\ (x\in X)$.
\end{prp}
\P By \ref{rdim+inf} we know that $\dim x=\sup (S_1)_x^{+\infty}\ (x\in X_1), \dim y=\sup(S_2)_y^{+\infty}\ (y\in X_2)$. The density of $S_1,S_2$ and the continuity of $g$ gives that $g(\sup P)=\sup(g(P))\ (P\subset S_1)$. But $g((S_1)_x^{+\infty})\subset(S_2)_{f(x)}^{+\infty}$ which gives the statement.
\Pe

We present a straightforward but important example for semi-isomorphim.

\begin{prp}Let ${\cal{D}}=\langle X,\mu_s:s\in S\rangle$ be a dimension structure. Set $\nu_s(x)=sign\big(\mu_s(x)\big)$ when $x\in X,s\in S$. Then ${\cal{D}}'=\langle X,\nu_s:s\in S\rangle$ is a dimension structure. Moreover ${\cal{D}}$ and ${\cal{D}}'$ are semi-isomorph for $f=id_X,g=id_S$. \Pes
\end{prp}

\subsection{On extensions}

The most problematic axiom of a dimension structure is (ax3) which requires the existence of an infimum. In this subsection we are going to investigate if a pre-structure satisfies (ax1) and (ax2) then how we can extend it to satisfies (ax3) as well.

\begin{df}Let ${\cal{D}}=\langle X,\mu_s:s\in S\rangle$ be given such that $X$ is a set, $\langle S,<\rangle$ is a partially ordered set, $\mu_s\ (s\in S)$ is a function $\mu_s:X\to\mathbb{R}^+\cup\{0,+\infty\}$ and axioms (ax1) and (ax2) hold. Then we call ${\cal{D}}$ a pre-dimension structure.
\end{df}

\begin{thm}\label{textdm}Let ${\cal{D}}=\langle X,\mu_s:s\in S\rangle$ be a pre-dimension structure with $S$ being a lattice. Then there exists an extension $\hat{\cal{D}}$ of ${\cal{D}}$ such that $\hat{\cal{D}}=\langle X,\mu_s:s\in \hat{S}\rangle$ is a dimension structure, $S\subset\hat{S}$ and if $p,q\in S$ then $p<_{S}q\iff p<_{\hat{S}}q$.
\end{thm}
\P If for $x\in X$ the infimum of $S_x$ does not exist then add a new element to $S$ i.e. let 
$$\hat{S}=S\cup\{\hat{x}:x\in X\text{ and }\not\exists\inf S_x\}$$
with the equivalence that $\hat{x}=\hat{y}$ if $S_x=S_y$ (or more precisely add the equivalence classes to $S$).

We have to define $\mu$ on the new elements. If $\hat{x}\in \hat{S}-S,\ y\in X$ then let
$$\mu_{\hat{x}}(y)=\sup\{\mu_s(y):s\in S_x\}.$$
Let us observe that this definition works for $p\in S$ too: $\mu_{p}(y)=\sup\{\mu_s(y):p\leq s\}$.

Let us define order on $\hat{S}$: If $t,u\in\hat{S}$ then let
$$t\leq_{\hat{S}} u\iff
\begin{cases}
\text{if }t,u\in S\text{ then }t\leq_{S}u\\
\text{if }t=\hat{x}\in\hat{S}-S,u\in S\text{ then }u\in S_x\\
\text{if }t\in S,u=\hat{y}\in\hat{S}-S\text{ then }p\in S_y\text{ implies that }t\leq_{S}p\\
\text{if }t=\hat{x},u=\hat{y}\in\hat{S}-S\text{ then }S_y\subset S_x.
\end{cases}$$
We have to check if it is an order. 

Reflexivity: We have to check points in $S$  and $\hat{S}-S$ and both cases are straightforward.

Antisymmetry: If $t,u\in S$ or $t,u\in\hat{S}-S$ then it is obvious. Let $\hat{x}\leq u,\ u\leq\hat{x}\ (u\in S)$. Then $u\in S_x$ and $\forall p\in S_x\ u\leq p$ which yields that $u=\inf S_x$ that is a contradiction. Let now $t\leq\hat{y},\ \hat{y}\leq t\ (t\in S)$. Then $\forall p\in S_y\ t\leq p$ and $t\in S_y$ which means that $t=\inf S_y$ -- a contradiction again.

Transitivity: Let assume that $t\leq u,u\leq w$. We have to check 8 cases. If $t,u,w\in S$ then it is clear that $t\leq w$ holds. If $t=\hat{x},u=\hat{y},w=\hat{z}\in\hat{S}-S$ then $S_z\subset S_y\subset S_x$ gives that $S_z\subset S_x$ which is $t\leq w$. If $t,u\in S,w=\hat{z}$ then $p\in S_z$ implies that $u\leq p$, and then $t\leq p$. If $t,w\in S,u=\hat{y}$ then $s\in S_y\Rightarrow t\leq s$ and $w\in S_y$ which yields that $t\leq w$. If $t=\hat{x},u,w\in S$ then $u\in S_x,u\leq w\Rightarrow w\in S_x$. If $t\in S,u=\hat{y},w=\hat{z}$ then $S_z\subset S_y$ and $s\in S_y\Rightarrow t\leq s$ hence $s\in S_z\Rightarrow t\leq s$. If $t=\hat{x},w=\hat{z},u\in S$ then $u\in S_x,\ s\in S_z\Rightarrow u\leq s$ that yields that $S_z\subset S_x$. If $t=\hat{x},u=\hat{y},w\in S$ then $S_y\subset S_x,\ w\in S_y$ which implies that $w\in S_x$.

\begin{itemize}
\item[(ax1)] Let $t<u$ and $\mu_t(x)<+\infty$. We have to check 4 cases. If $t,u\in S$ then it is (ax1) in ${\cal{D}}$. If $t=\hat{x},u\in S$ then $u\in S_x$ i.e. $\mu_u(x)=0$. If $t\in S,u=\hat{y}$ then $s\in S_y\Rightarrow t< s$ hence $\mu_s(x)=0$ for all $s\in S_y$ which gives that $\mu_{\hat{y}}(x)=0$. If $t=\hat{x},u=\hat{y}$ then $S_y\subset S_x$ hence $\mu_s(x)=0$ for all $s\in S_y$ which gives that $\mu_{\hat{y}}(x)=0$.
\item[(ax2)] Let $0<\mu_t(y)<+\infty,\mu_u(y)<+\infty$. We have to check 4 cases. If $t,u\in S$ then it is (ax2) in ${\cal{D}}$. If $t=\hat{x},u\in S$ then $\forall p\in S_x\ \mu_p(y)<+\infty$ and $\exists q\in S_x$ such that $0<\mu_q(y)<+\infty$ which gives that $q\leq u$ hence $u\in S_x$ i.e. $\hat{x}\leq u$. 
If $t\in S,u=\hat{x}$ then $\forall p\in S_x\ \mu_p(x)<+\infty$ which means that each $p\in S_x$ is comparable with $t$ therefore $t\leq p$  i.e. $t\leq\hat{x}$. If $t=\hat{x},u=\hat{z}$ then $\forall p\in S_x\ \mu_p(x)<+\infty$ and $\exists q\in S_x$ such that $\mu_q(y)>0$ and $\forall r\in S_z\ \mu_r(y)<+\infty$ which gives that $\forall r\in S_z\ q\leq r$ hence $S_z\subset S_x$ i.e. $\hat{x}\leq\hat{z}$.
\item[(ax3)] Let $$\hat{S}_x=\{s\in\hat{S}:\mu_s(x)<+\infty\}=S_x\cup\{\hat{y}\in \hat{S}-S:\mu_{\hat{y}}(x)<+\infty\}.$$
Let $\inf_S,\ \inf_{\hat{S}}$ denote the infimum in $S$ and $\hat{S}$ respectively.

We have two cases.

(1) If $\inf_S S_x$ does not exist in $S$. Now we show that $\inf_{\hat{S}} \hat{S}_x=\hat{x}$. If $p\in S_x$ then $\hat{x}\leq p$ by definition of $\leq_{\hat{S}}$. If $\mu_{\hat{y}}(x)<+\infty$ then $\sup\{\mu_s(x):s\in S_y\}<+\infty$ which gives that $\forall s\in S_y\ \mu_s(x)<+\infty$ hence $s\in S_x$ i.e. $S_y\subset S_x$ that is $\hat{x}\leq\hat{y}$.

If $t\in S,\forall s\in\hat{S}_x\ t\leq_{\hat{S}} s$ then $\forall s\in S_x\ t\leq s$ i.e. $t\leq \hat{x}$. 
If $\hat{y}\in \hat{S}-S,\forall s\in\hat{S}_x\ \hat{y}\leq s$ then $\forall s\in S_x\ \hat{y}\leq s$ hence $s\in S_y$ which gives that $S_x\subset S_y$ i.e. $\hat{y}\leq \hat{x}$.

(2) If $s_0=\inf_S S_x$ exists in $S$. In this case we show that $s_0=\inf_{\hat{S}}\hat{S}_x$. If $s\in S_x$ then $s_0\leq s$. If $\hat{y}\in\hat{S}_x-S$ then $\mu_{\hat{y}}(x)<+\infty$ i.e. $\forall p\in S_y\ \mu_p(x)<+\infty$ hence $p\in S_x$ which gives that $s_0\leq p$ and we get that $s_0\leq \hat{y}$.

If $s\in S$ and $\forall t\in\hat{S}_x\ s\leq t$ then $\forall t\in S_x\ s\leq t$. Therefore $s\leq s_0$. If $\forall t\in\hat{S}_x\ \hat{y}\leq t$ then $\forall t\in S_x\ \hat{y}\leq t$ that is $t\in S_y$ that gives that $S_x\subset S_y$. Then $S_x=S_y$ would be a contradiction hence there exists $s\in S_y-S_x$. From lattice theory we get that $\inf_{S} S_x\cup\{s\}=s_0\wedge s$. From part (1) $\hat{y}=\inf_{\hat{S}}\hat{S}_y$ that gives that $\hat{y}\leq s_0\wedge s\leq s_0$.
\end{itemize}
\Pes

\begin{prp}Let ${\cal{D}}=\langle X,\mu_s:s\in S\rangle$ be a pre-dimension structure and $\hat{\cal{D}}$ be its extension as described in \ref{textdm}. If ${\cal{D}}$ is principal then so is $\hat{\cal{D}}$.
\end{prp}
\P We use the notations introduced in the proof of \ref{textdm}.

There are two cases we have to manage. During that we will apply \ref{pprinceq}.

(1) If $s_0=\inf_{\hat{S}}\hat{S}_x\in S$. If $s_0<s\in S$ then $\mu_s(x)=0$. If $s_0<\hat{y}$ then $\forall p\in S_y\ s_0<p$ hence $\mu_p(x)=0$ which gives that $\mu_{\hat{y}}(x)=0$.

(2) If $\hat{x}=\inf_{\hat{S}}\hat{S}_x\in \hat{S}-S$. If $\hat{x}<s\ (s\in S)$ then $s\in S_x$ that is $\mu_s(x)<+\infty$ but $s\ne\inf_S S_x$ hence $\mu_s(x)=0$. If $\hat{x}<\hat{y}$ then $S_y\subset S_x$ which gives that $\forall p\in S_y\ \mu_p(x)<+\infty$ but $p\ne\inf_S S_x$ therefore $\mu_p(x)=0\ (\forall p\in S_y)$ which implies that $\mu_{\hat{y}}(x)=0$.
\Pe

Now we investigate how unique the previous extension is. 

\begin{thm}Let ${\cal{D}}=\langle X,\mu_s:s\in S\rangle$ be a pre-dimension structure with $S$ being a lattice. Let ${\cal{D}}$ be principal. Let $\hat{\cal{D}}$ be its extension as described in \ref{textdm}. Let ${\cal{D}}'=\langle X,\mu_s:s\in S'\rangle$ be a dimension structure such that ${\cal{D}}$ is one of its substructures. Then there is $f:\hat{S}\to S'$ that preserves order, injective and keeps $S$ pointwise fixed ($s\in S\subset \hat{S}\Rightarrow f(s)=s\in S'$). Moreover $\mu^{{\cal{D}}'}_s(x)=\mu^{\hat{\cal{D}}}_s(x)$ if $s\in S$, and $\mu^{{\cal{D}}'}_{f(\hat{x})}(x)\geq\sup\{\mu^{\hat{\cal{D}}}_s(x):s\in S_x\}$.
\end{thm}
\P Let $$f(t)=
\begin{cases}
t&\text{if }t\in S\\
\inf_{S'}S_x&\text{if }t\in\hat{S}-S,\ t=\hat{x}.
\end{cases}$$
We show that $f$ preserves order. Let $t,u\in\hat{S}\ t<u$. If $t,u\in S$ then it is clear. If $t=\hat{x},u\in S$ then $u\in S_x$ hence $\inf_{S'}S_x\leq u$. If $t\in S,u=\hat{y}$ then $\forall s\in S_y\ t\leq s$ therefore $t\leq \inf_{S'}S_y$. If $t=\hat{x},u=\hat{y}$ then $S_y\subset S_x$ which yields that $\inf_{S'}S_x\leq\inf_{S'}S_y$.

Now we prove that $f$ is injective. Let $t,u\in\hat{S},\ t\ne u$. If $t,u\in S$ then it is obvious. If $t=\hat{x},u\in S$ then $\inf_{S'}S_x\ne u$ because let us assume that they were equal. Then $\inf_{S'}S_x\in S$ would hold which would give that $t=u$ -- a contradiction. If $t=\hat{x},u=\hat{y}$ then assume that $\inf_{S'}S_x=\inf_{S'}S_y$. By principality it would yield that $S_x=S_y$ hence $\hat{x}=\hat{y}$ would hold -- a contradiction.

Evidently $\mu^{{\cal{D}}'}_s(x)=\mu^{\hat{\cal{D}}}_s(x)\ (s\in S)$ holds because ${\cal{D}}$ is a substructure of both $\hat{\cal{D}}$ and ${\cal{D}}'$. If $s\in S_x$ then $\mu^{{\cal{D}}'}_{f(\hat{x})}(x)\geq\mu^{{\cal{D}}'}_s(x)$ but $\mu^{{\cal{D}}'}_s(x)=\mu^{\hat{\cal{D}}}_s(x)$ which gives the last statement.
\Pes


{\footnotesize

\noindent

\noindent E-mail: alosonczi1@gmail.com\\


\begin{thebibliography}{www}

\bibitem{billings} P. Billingsley, {\em Probability and Measure}, John Wiley \& Sons, vol. 939  (2012).

\bibitem{erd} R. Engelking, {\em Dimension theory}, North-Holland Publishing Company (1978).

\bibitem{ert} R. Engelking, {\em General topology}, Sigma Series in Pure Mathematics, vol. 6. (1989).

\bibitem{fk} K. Falconer, {\em Fractal geometry: mathematical foundations and applications}, John Wiley \& Sons (2004).

\bibitem{hf} F. Hausdorff, {\em Dimension and outer measure}, Mathematische Annalen {\em 79.1} (1918), 157–-179.

\bibitem{lambm} A. Losonczi, {\em Measures by means, means by measures}, arXiv prepint

\bibitem{mp} B. B. Mandelbrot, R. Pignoni, {\em The fractal geometry of nature}, New York: WH freeman Vol. 173 (1983).

\bibitem{rs} S. Roman, {\em Lattices and ordered sets}, Springer Science \& Business Media (2008).


\end{thebibliography}
\end{document}